\documentclass{article}

\usepackage{amsmath,amssymb,amsfonts,amsthm,amscd}

\newtheorem{prop}{Proposition}

\title{A new invariant under congruence\\ of nonsingular matrices}

\author{Kiyoshi Shirayanagi and Yuji Kobayashi}

\date{}
  
\begin{document}
\maketitle

\footnote[0]{\hspace{-4ex}{\it Keywords}:  Linear algebra, Invariant of matrices, Congruence, Zeropotent algebras}
\footnote[0]{\hspace{-4ex}{\it 2010 Mathematics Subject Classification}: Primary 15A15; Secondary 15A72}

\footnote[0]{\hspace{-4ex}Research of the first author is supported in part by JSPS KAKENHI Grant Number JP18K11172.}
\footnote[0]{\hspace{-4ex}Research of the second author was supported in part by JSPS KAKENHI Grant Number JP25400120.}

\begin{abstract}
For a nonsingular matrix $A$, we propose the form $Tr(^t\!A A^{-1})$, the trace of the product of its transpose and inverse, 
as a new invariant under congruence of nonsingular matrices.

\end{abstract}

\section{Introduction}
Let $K$ be a field. We denote the set of all $n\times n$ matrices over $K$ by $M(n,K)$,
and the set of all $n\times n$ nonsingular matrices over $K$ by $GL(n,K)$.
For $A, B\in M(n,K)$, $A$ is {\it similar} to $B$ if there exists $P\in GL(n,K)$ such that
$P^{-1}AP=B$, and $A$ is {\it congruent} to $B$ if there exists $P\in GL(n,K)$ such that
$^t\!PAP=B$, where $^t\!P$ is the transpose of $P$. A map $f: M(n,K) \rightarrow K$ is an {\it invariant under similarity of matrices} if for any $A\in M(n,K)$,
$f(P^{-1}AP)=f(A)$ holds for all $P\in GL(n,K)$.
Also, a map $g: M(n,K) \rightarrow K$ is an {\it invariant under congruence of matrices} if for any $A\in M(n,K)$,
$g(^t\!PAP)=g(A)$ holds for all $P\in GL(n,K)$, and 
$h: GL(n,K) \rightarrow K$ is an {\it invariant under congruence of nonsingular matrices} if for any $A\in GL(n,K)$,
$h(^t\!PAP)=h(A)$ holds for all $P\in GL(n,K)$.

As is well known, there are many invariants under similarity of matrices, such as trace, determinant and other coefficients of 
the minimal polynomial. In case the characteristic is zero, rank is also an example.
Moreover, much research on similarity (or {\it simultaneous conjugation}) has been done from the viewpoint of invariant theory for a group action on tuples of matrices, for example, see \cite{procesi}.
On the other hand, few invariants under congruence of matrices or nonsingular matrices are known except for rank.
Restricted to symmetric matrices over the reals, another exception would be from Sylvester's law of inertia, which states that the numbers of positive, negative, and zero eigenvalues
are all invariants under congruence of real symmetric matrices. 
In this note, we propose a new invariant under congruence of general nonsingular matrices.
The form is $\sigma(A)=Tr(^t\!A A^{-1})$, where $Tr$ denotes trace.
\section{Invariance and other properties of $\sigma$}



First of all, let us show that $\sigma$ is an invariant under congruence of nonsingular matrices.
In fact, for any nonsingular matrix $X$,

\begin{eqnarray*}
\sigma(^t\!XAX)&=&Tr(^t(^t\!XAX) (^t\!XAX)^{-1})\\
&=&Tr((^t\!X\,^t\!AX)(X^{-1}A^{-1}(^t\!X)^{-1}))\\
&=&Tr(^t\!X(^t\!AA^{-1})(^t\!X)^{-1})\\
&=&\sigma(A).
\end{eqnarray*}

As an immediate corollary, any polynomial in $\sigma(A)$ over $K$ is also an invariant under congruence of nonsingular matrices. 
Moreover, any rational function in $\sigma(A)$ over $K$ is an invariant 
under congruence of nonsingular matrices for which its denominator does not vanish.


Next, let us describe some other properties of $\sigma$.

\begin{prop}\label{prop:properties}
\noindent For any $n\times n$ nonsingular matrix $A$ over $K$, the following equalities hold.
\begin{enumerate}
\renewcommand{\labelenumi}{\rm{(\arabic{enumi})}}
\item\label{enu:symmetry} $\sigma(^t\!A)=\sigma(A)$.
\item\label{enu:inverse} $\sigma(A^{-1})=\sigma(A)$.
\item\label{enu:adj} $\sigma(adj(A))=\sigma(A)$, where $adj(A)$ is the adjugate of $A$.
\item\label{enu:scalar} $\sigma(cA)=\sigma(A)$ for any $c\in K\setminus\{0\}$.
\item\label{enu:symmetric} $\sigma(A)=n$ if $A$ is symmetric.
\end{enumerate}
\end{prop}

\noindent{\it Proof.}

(\ref{enu:symmetry})
\begin{eqnarray*}
\sigma(^t\!A)&=&Tr(^t(^t\!A)\, (^t\!A)^{-1})\\
&=&Tr(A\,^t(A^{-1}))\\
&=&Tr(A^{-1}(^t\!A))\\
&=&\sigma(A)
\end{eqnarray*}

(\ref{enu:inverse})
\begin{eqnarray*}
\sigma(A^{-1})&=&Tr(^t(A^{-1})\,(A^{-1})^{-1})\\
&=&Tr(^t(A^{-1})\,A)\\
&=&\sigma(A)
\end{eqnarray*}

(\ref{enu:adj})
\begin{eqnarray*}
\sigma(adj(A))&=&Tr(^t(adj(A))\,(adj(A))^{-1})\\
&=&Tr(adj(^t\!A)\,(adj(A))^{-1})\\
&=&Tr((|^t\!A|(^t\!A)^{-1})\,(|A|^{-1}A))\\
&=&Tr((^t\!A)^{-1}\,A)\\
&=&Tr(^t(A^{-1})\,A)\\
&=&\sigma(A)
\end{eqnarray*}

(\ref{enu:scalar}) For any $c\in K\setminus\{0\}$,
\begin{eqnarray*}
\sigma(cA)&=&Tr(^t(cA)\,(cA)^{-1})\\
&=&Tr(c\,^t\!A\,c^{-1}A^{-1})\\
&=&\sigma(A)
\end{eqnarray*}

(\ref{enu:symmetric}) If $^t\!A=A$, then
\begin{eqnarray*}
\sigma(A)&=&Tr(AA^{-1})\\
&=&Tr(E)\\
&=&n
\end{eqnarray*}

\section{Background on $\sigma$}
Let us briefly describe the background that led to the discovery of $\sigma$. The authors of the present note with Sin-Ei Takahasi and Makoto Tsukada
tried to classify three-dimensional {\it zeropotent} algebras up to isomorphism, where a zeropotent algebra is a nonassociative algebra 
in which the square of any element is zero, see \cite{KSTT,STTK}. We there expressed an algebra by its structure constants. For a three-dimensional zeropotent algebra $A$,
let $\{e_1, e_2, e_3\}$ be a linear base of $A$. By zeropotency of $A$ it suffices to consider the structure constants for $e_2e_3, e_3e_1, e_1e_2$,
namely, the algebra can be identified with the $3\times 3$ matrix $A$ such that $^t(e_2e_3\,\, e_3e_1\,\, e_1e_2)=A\,^t(e_1\,\, e_2\,\, e_3)$. We hereafter use the same symbol $A$ both for the matrix and for the algebra.

The following proposition gives a criterion for two zeropotent algebras to be isomorphic.
\begin{prop}[\cite{KSTT}]\label{prop:equiv}
Let $A$ and $A'$ be three-dimensional zeropotent algebras over $K$. Then, 
$A$ and $A'$ are isomorphic if and only if there is a nonsingular matrix $X$ satisfying 
\begin{equation*}\label{eq:equiv}
A'=\frac{1}{|X|}\,\,  \\^t\!XAX.
\end{equation*}
In particular, 
if $K$ is algebraically closed, then $A$ and $A'$ are isomorphic if and only if there is a nonsingular matrix $X$ satisfying $A'=\,^t\!XAX$, that is, $A$ and $A'$ are congruent.
\end{prop}

Our strategy of classification was to divide zeropotent algebras into {\it curly} algebras and {\it straight} algebras. A zeropotent algebra is curly if 
the product of any two elements lies in the space spanned by these elements, otherwise straight.
Let $K$ be an algebraically closed field.
A straight algebra of rank 3 can be expressed by the {\it canonical form} 
\[
A(a, b, c)=\begin{pmatrix}1&a&b\\0&1&c\\0&0&1\end{pmatrix}
\]
with $a, b, c\in K$. 

During the course of classifying the straight algebras, we encountered the quantity $D=a^2+b^2+c^2-abc$, which seems strange 
because it is not homogeneous with respect to $a$, $b$, and $c$.
That is,
in case where $D\ne 0$ and $a^2+b^2-abc\ne 0$, we showed that $A(a,b,c)$ is congruent to $A(\sqrt{D},0,0)$.
From this, it turned out that $D$ is an invariant under congruence of some $A(a,b,c)$-type matrices.
Inspired by this fact, we got interested in elucidating the origin of $D$.
From $D$, later, by backward reasoning  from the canonical form to a general matrix using computer, 
the second author derived a more general homogeneous quantity 
$\dfrac{1}{|A|}\{
a_{13}^2a_{22}+3a_{12}a_{23}a_{31}-a_{21}a_{23}a_{31}+a_{22}a_{31}^2+a_{11}(a_{23}-a_{32})^2-a_{12}a_{31}a_{32}-a_{21}a_{31}a_{32}-a_{13}(a_{21}a_{23}+2a_{22}a_{31}-3a_{21}a_{32}+a_{12}(a_{23}+a_{32}))+a_{12}^2a_{33}-2a_{12}a_{21}a_{33}+a_{21}^2a_{33}\}$,
for a $3\times 3$ matrix $A=(a_{ij})$ of rank 3. Let us denote this quantity by $\kappa(A)$.
By computer we confirmed that $\kappa$ is certainly an invariant under congruence of general nonsingular matrices, but could not prove it mathematically.
Since the original expression of $\kappa(A)$ was awkward and complicated,
by an enormous amount of transformations of trial and error, the first author obtained the more beautiful form
$3-\dfrac{1}{|A|}\displaystyle\sum_{i,j}(a_{ij}\widetilde{a_{ji}})$, 
where $\widetilde{a_{ij}}$ denotes the $(i,j)$ cofactor, that is, the product of $(-1)^{i+j}$ and the $(i,j)$ minor.
Denoting $\dfrac{1}{|A|}\displaystyle\sum_{i,j}(a_{ij}\widetilde{a_{ji}})$ by $\sigma(A)$,
by further reasoning he had $\sigma(A)=\dfrac{1}{|A|}Tr(A\, adj(^t\!A))$, and finally by using  $^t\!A(adj(^t\!A))=|^t\!A|E$, 
reached the simple form $Tr(^t\!A A^{-1})$ as initially defined
and then $\kappa(A)=3-Tr(^t\!A A^{-1})$. Through this expression we can now mathematically prove that
$\kappa$ is an invariant under congruence of nonsingular matrices.

Beyond the topic of three-dimensional zeropotent algebras, we found that in general, for any integer $n\geq 1$,
$\sigma$ is an invariant under congruence of $n\times n$ nonsingular matrices over a general field.



\bibliographystyle{amsplain}

\begin{thebibliography}{[1]}
\bibitem{KSTT} Y. Kobayashi, K. Shirayanagi, S.-E. Takahasi and M. Tsukada, Classification of three-dimensional zeropotent algebras over an algebraically closed field, {\it Comm. Algebra}, {\bf 45} 12 (2017), 5037--5052.
\bibitem{procesi} C. Procesi, The invariant theory of $n\times n$ matrices, {\it Adv. Math.}, {\bf 19} (1976), 306--381.
\bibitem{STTK} K. Shirayanagi, S.-E. Takahasi, M. Tsukada and Y. Kobayashi, Classification of three-dimensional zeropotent algebras over the real number field, {\it Comm. Algebra}, {\bf 46} 11 (2018), 4663--4681.

\end{thebibliography}

\noindent
K. Shirayanagi and Y. Kobayashi\\
Department of Information Science, Toho University, Miyama 2-2-1 Funabashi, Chiba 274-8510, Japan\\ kiyoshi.shirayanagi@is.sci.toho-u.ac.jp and kobayasi@is.sci.toho-u.ac.jp

\end{document}